# Optimal Pricing Strategies for Heterogeneous Customers in Dual-Channel Closed-Loop Supply Chains: A Modeling Approach

Yang XIAO [1*], Hisashi KURATA [2], Ting WANG [2]

[1] Faculty of Economics, Saitama University, Japan

[2] Graduate School of Internation Social Sciences, Yokohama National University, Japan

* Corresponding author: xiaoyangq1996@gmail.com

**Abstract.** Dual-channel closed-loop supply chains (DCCLSCs) play a vital role in attaining both sustainability and profitability. This paper introduces a game-theoretic model to analyze optimal pricing strategies for primary and replacement customers within three distinct recycling frameworks: manufacturer-led, retailer-led, and collaborative recycling. The model identifies equilibrium pricing and subsidy decisions for each scenario, considering the primary customer's preference for the direct channel and the specific roles in recycling. The findings indicate that manufacturers tend to set lower prices in direct channels compared to retailers, aiming to stimulate demand and promote trade-ins. Manufacturer-led recycling initiatives result in stable pricing, whereas retailer-led recycling necessitates higher subsidies. Collaborative recycling strategies yield lower prices and an increase in trade-ins. Primary customers' preference for the direct channel significantly impacts pricing strategies, with a stronger preference leading to lower direct-channel prices and higher manufacturer subsidies. This paper contributes to the field by incorporating primary customer channel preferences and diverse recycling frameworks into DCCLSC pricing models. These insights assist manufacturers and retailers in adjusting pricing strategies and trade-in incentives according to primary customer preferences and associated costs, thereby enhancing profitability and recycling efficiency within DCCLSCs.

## 1  Introduction

The field of supply chain management has significantly evolved due to the intersection of digital transformation and the increasing need for environmental sustainability. A notable development within this context is the emergence of dual-channel closed-loop supply chains (DCCLSCs), defined as supply chain systems that integrate both online and physical sales channels alongside forward and reverse logistics processes. DCCLSCs have become essential business models due to their adaptability to consumer preferences and environmental considerations. While traditional studies have extensively explored pricing strategies within dual-channel contexts (Huang et al., 2009), these analyses predominantly focused on homogeneous customer groups. However, contemporary markets emphasize customer heterogeneity, differentiating between primary customers—those who purchase new products without returning used ones—and replacement customers, who engage in trade-in programs by returning used items for discounts on new purchases (Zhou et al., 2021).

Closed-loop supply chains (CLSCs), which incorporate product returns, recycling, and remanufacturing processes, have gained increasing prominence due to their dual benefit of sustainability and profitability (Daniel et al., 2002; Guide et al., 2003; Srivastava, 2007). Unlike traditional supply chains, CLSCs emphasize a circular economy through the integration of forward and reverse logistics flows, significantly reducing environmental impact. Industry examples include Unilever's sustainable sourcing practices and Patagonia's repair-and-reuse initiatives, underscoring the economic and environmental viability of CLSC practices (Giva, 2022; Greenfield, 2023). Recent academic contributions further highlight methodologies such as multi-criteria decision-making, fuzzy modeling, blockchain technology, and robust optimization frameworks to address complexities within CLSCs, thereby enhancing operational efficiency and managing uncertainties (Ali et al., 2020; Ali et al., 2021; Goli et al., 2023; Lotfi et al., 2023).

Despite these advancements, comprehensive analyses of optimal pricing strategies in DCCLSCs considering heterogeneous customer segments and varied recycling frameworks remain limited. Although elements such as dual-channel structures (Chiang et al., 2005; Arya et al., 2008), recycling processes (Savaskan et al., 2004; Ferrer & Swaminathan, 2006), and customer heterogeneity (Chen et al., 2017; Overby & Lee, 2006) have been individually explored, few studies have examined their combined effects thoroughly. Recent investigations have started addressing some aspects: Wu et al. (2025) analyzed green innovation and big data marketing within CLSCs, demonstrating enhanced profitability through coordinated strategies. Similarly, Chen and Wu (2024) studied offline preferences and service levels in pricing decisions, highlighting the importance of these variables. Further research by Yu et al. (2024, 2025) emphasized dynamic pricing and recycling decisions, while He et al. (2024) and Wei (2024) respectively focused on competitive dynamics among collectors and governmental subsidy impacts. Additional contributions examined customer segmentation (Rahmani & Pashapour, 2024), psychological influences on recycling behavior (Huang et al., 2024), effectiveness of dual collection channels and wholesale price discounts (Pal et al., 2024), bundling strategies and advertising impacts (Hadadi et al., 2024), and the implications of asymmetric information and uncertainty (Beranek & Buscher, 2024).

However, the intricate interactions among customer heterogeneity—particularly distinguishing primary and replacement customers—channel preferences, and diverse recycling structures (manufacturer-led, retailer-led, and collaborative recycling) have not been integrated comprehensively in previous research. This paper addresses this critical gap by systematically analyzing how these combined dimensions influence optimal pricing and recycling subsidy strategies within a DCCLSC framework. **Table 1** summarizes this paper's unique contributions compared to previous related literatures.

Specifically, this paper investigates three core research questions:

(1) How do pricing strategies differ between primary and replacement customers in DCCLSCs?

(2) How do manufacturer-led, retailer-led, and collaborative recycling frameworks impact these pricing strategies?

(3) How do interactions between customer segments, channel preferences, and recycling structures affect sustainability objectives?

This paper develops a Manufacturer-Stackelberg game-theoretic model (Choi, 1991) to analyze optimal pricing strategies within DCCLSCs. The model explicitly incorporates differentiated demand functions reflecting the valuation and purchase behaviors of primary and replacement customers across direct and indirect sales channels. By solving for equilibrium conditions under three distinct recycling frameworks—manufacturer-led, retailer-led, and collaborative recycling—the analysis systematically identifies optimal prices and subsidies. Utilizing this modeling approach, this paper finds that manufacturers typically set lower direct channel prices compared to retailers to attract replacement customers and stimulate trade-in activities. Additionally, the analysis reveals a significant sensitivity of pricing strategies to variations in customer heterogeneity and recycling frameworks, demonstrating that collaborative recycling structures generally yield lower prices, higher trade-in rates, and enhanced recycling efficiency. Importantly, the results indicate that manufacturer-led recycling frameworks provide stable pricing environments, whereas retailer-led frameworks require higher subsidies to incentivize trade-ins. These findings offer critical managerial insights into effectively balancing economic performance and sustainability objectives within DCCLSCs.

The remainder of this paper is structured as follows: Section 2 outlines the model's assumptions and key variables. Section 3 describes the methodology and analytical framework employed for determining pricing strategies. Numerical analyses in Section 4 illustrate practical applications and comparative insights. Section 5 concludes by summarizing key findings, managerial implications, and identifying avenues for future research.

**TABLE 1.** Comparative summary of related literatures and this paper

| Study | Dual-channel Pricing | Customer Heterogeneity | Channel Preference | Recycling Frameworks | Subsidy and Incentive Analysis | Dynamic/Static Model | Focus of Analysis |
|---|---|---|---|---|---|---|---|
| Wu et al. (2025) | ✓ | ✗ | ✗ | Single (Centralized) | ✓ (Cost-sharing) | Dynamic (Differential Game) | Green innovation, big data |
| Chen & Wu (2024) | ✓ | ✗ | ✓ (Offline) | ✗ | ✗ | Static (Stackelberg) | Offline preference, service levels |
| Yu et al. (2024) | ✓ | ✗ | ✗ | Dual-channel (Dynamic) | ✓ | Dynamic | Pricing and recycling coordination |
| Yu et al. (2025) | ✓ | ✗ | ✗ | Platform-based modes | ✓ | Static | Trade-in operation decisions |
| He et al. (2024) | ✓ | ✗ | ✓ (Convenience) | Multiple collectors | ✓ | Static | Collection competition |
| Wei (2024) | ✗ | ✓ (Remanufactured acceptance) | ✗ | Single-channel | ✓ (Double subsidies) | Static | Subsidy effects on remanufacturing |
| Rahmani & Pashapour (2024) | ✓ | ✓ (Customer segments) | ✗ | Dual-channel | ✗ | Static (optimization) | Customer segmentation |
| Huang et al. (2024) | ✗ | ✓ (Psychological factors) | ✗ | Single-channel | ✓ (Social emotions) | Dynamic (Evolutionary) | Social influence on recycling |
| Pal et al. (2024) | ✓ | ✗ | ✗ | Dual collection | ✓ (Price discount) | Static | Green collection effort |
| Hadadi et al. (2024) | ✓ | ✗ | ✗ | Dual-channel bundling | ✓ (Joint advertising) | Static | Bundling and advertising |
| Beranek & Buscher (2024) | ✓ | ✗ | ✗ | Single collection (Asymmetric info.) | ✓ (Info. asymmetry) | Dynamic (Two-period game) | Asymmetric information |
| This Paper | ✓ | ✓ (Primary vs. Replacement) | ✓ (Direct channel preference) | Multiple (Manufacturer-led, Retailer-led, Collaborative) | ✓ (Trade-in subsidies and transfer pricing) | Static (Stackelberg equilibrium) | Integrated effects of customer heterogeneity, channel preference, and multiple recycling frameworks |

## 2 Model Description

### 2.1 Problem Description

This paper investigates a DCCLSC model consisting of a single manufacturer and a single retailer, as illustrated in **FIGURE 1**. This model facilitates collaboration between the manufacturer and retailer in the domains of new-product distribution and trade-in services. It effectively represents a typical dual-channel scenario, allowing for the exploration of various recycling frameworks under differing cost parameters and customer preferences. The DCCLSC is analyzed through three distinct recycling models: Manufacturer-led recycling (Model M), Retailer-led recycling (Model R), and Co-recycling (Model MR). The supply chain serves two customer segments: primary customers, who directly purchase new items directly, and replacement customers, who trade in used items for remanufacturing and then purchase new items.

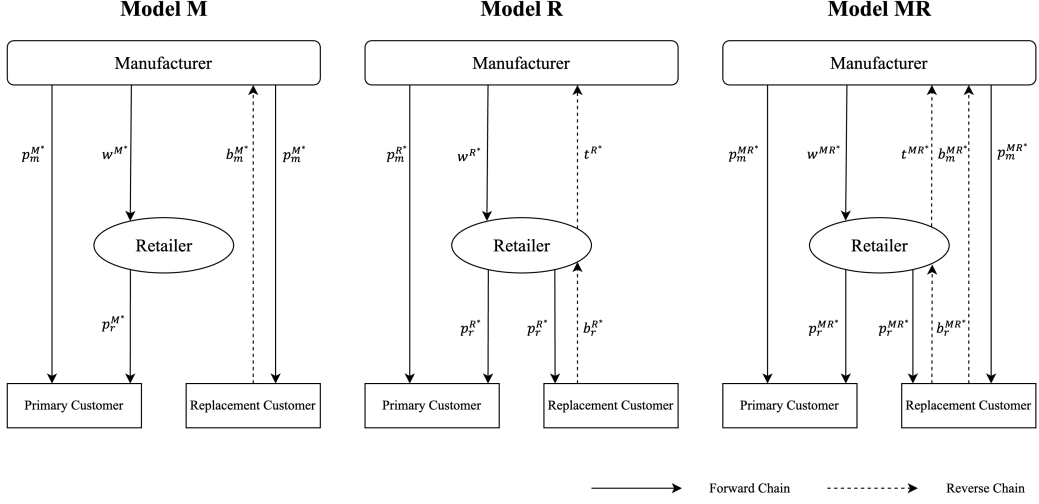

**FIGURE 1.** Three Recycling Models

In the forward chain, the manufacturer supplies new items to the retailer at a unit wholesale price of $w^{M^*}$, $w^{R^*}$, or $w^{MR^*}$, as detailed in **TABLE 2**. Both parties then sell the same new items through two channels: the manufacturer's direct channel at a unit selling price of $p_m^{M^*}$, $p_m^{R^*}$, or $p_m^{MR^*}$, and the retailer's indirect channel at a unit selling price of $p_r^{M^*}$, $p_r^{R^*}$, or $p_r^{MR^*}$. In the reverse chain, three distinct models are proposed for collecting used items from primary customers. In the Model M, the manufacturer directly collects used items by providing a unit subsidy denoted as $b_m^{M^*}$. In the Model R, the retailer is responsible for collecting used items by offering a unit subsidy $b_r^{R^*}$, while the manufacturer reimburses the retailer's collection efforts through a unit transfer price $t^{R^*}$, where $t^{R^*} \geq b_r^{R^*}$. In the Model MR, both the manufacturer and the retailer independently collect used items, offering subsidies $b_m^{MR^*}$ and $b_r^{MR^*}$, respectively. Under this model, the manufacturer also pays a unit transfer price $t^{MR^*}$, such that $t^{MR^*} \geq b_r^{MR^*}$.

**TABLE 2.** Decision Variables

| | |
|---|---|
| $p_m^i$, $p_r^i$ | The unit selling price for new items in the direct and indirect channels, $i \in \{M, R, MR\}$ |
| $b_m^i$, $b_r^i$ | The unit subsidy specified by the manufacturer and retailer for used items, $i \in \{M, R, MR\}$ |
| $w^i$ | The unit wholesale price for new items, $i \in \{M, R, MR\}$ |
| $t^i$ | The unit transfer price specified by the manufacturer for used items, $i \in \{M, R, MR\}$ |

## 2.2 Model Assumptions

The DCCLSC model under examination encompasses a multitude of decision variables, customer segments, and reverse logistics processes, which may become analytically cumbersome if excessively intricate. Consequently, the subsequent assumptions are employed to limit the scope to fundamental components, thereby ensuring the manageability of cost structures, strategic interactions, and customer behaviors. Each assumption is consistent with empirical observations in dual-channel supply chains, where the advantages of remanufacturing costs, the differentiation of customer types, and interfirm coordination are pivotal. TABLE 3 comprehensively lists and defines the key parameters employed in our analysis. The utilities ($U_j^i$) and market sizes ($q_j^i$) for customers who purchase new items are explicitly segmented across different recycling scenarios ($j \in \{1, 2, 3, 4\}$) and customer-channel combinations. Specifically, the index set $j$ clearly captures distinct purchasing behaviors and preferences: $j = 1$: Represents primary customers purchasing new items directly through the manufacturer's channel; $j = 2$: Represents primary customers purchasing new items indirectly through the retailer's channel; $j = 3$: Represents replacement customers purchasing new items via direct channel by trading in used items under the manufacturer's subsidy in Models M and MR, and via the indirect channel under the retailer's subsidy in Model R; $j = 4$: Represents replacement customers purchasing new items via indirect channel by trading in used items under the retailer's subsidy in Model MR.

By explicitly differentiating these four market segments, our analysis accurately captures the diverse purchasing behaviors and preferences of customers, thus ensuring a realistic representation of market dynamics. The profit functions ($\Pi_m^i$, $\Pi_r^i$, $\Pi_s^i$) for the manufacturer, retailer, and the entire supply chain, respectively, are precisely defined to facilitate accurate economic evaluations across each recycling scenario. Moreover, the parameter $\alpha$, representing primary customers' preference for the direct sales channel ($0 < \alpha < 1$), explicitly models observed customer heterogeneity in contemporary dual-channel markets. The valuation coefficients $v$ (primary customers) and $u$ (replacement customers) reflect the distinct perceived values for new products, allowing nuanced differentiation in customer segmentation and associated strategic implications. Additionally, parameters such as the government subsidy $s$, production costs ($c_m$, $c_r$), and unit cost savings from remanufacturing ($\Delta = c_m - c_r$) are grounded in realistic industrial observations and supported by established literature (Savaskan et al., 2004; Wei, 2024), ensuring the practical relevance and economic validity of our modeling assumptions.

**Assumption 1.** *Under a trade-in policy, customers are categorized into two groups: "primary" customers, who purchase new items without returning used ones, and "replacement" customers, who return used items while purchasing new ones.*

**Assumption 2.** *A manufacturer-led Stackelberg game framework is applied across all three models $i \in \{M, R, MR\}$. In this framework, the manufacturer (leader) first determines the variables $\{w^{i^*}, p_m^{i^*}, b_m^{i^*}, t^{i^*}\}$, including the wholesale price, direct selling price, transfer price, and manufacturer's subsidy. Subsequently, the retailer (follower) determines the variables $\{p_r^{i^*}, b_r^{i^*}\}$, which include the indirect selling price and retailer's subsidy, as illustrated in Figure 2.*

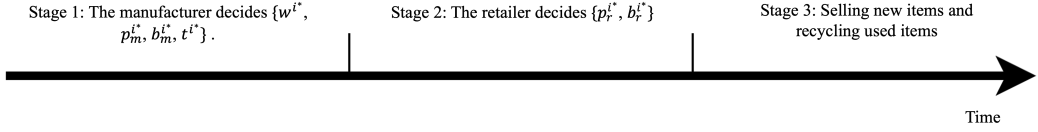

**FIGURE 2.** The Decision Sequence

**TABLE 3.** Definition of Parameters

| | |
|---|---|
| $U_j^i$ | The utility for customers who purchase a new item, $i \in \{M, R, MR\}$, $j \in \{1, 2, 3, 4\}$ |
| $q_j^i$ | The market size, $i \in \{M, R, MR\}$, $j \in \{1, 2, 3, 4\}$ |
| $\Pi_m^i, \Pi_r^i, \Pi_s^i$ | Profit for the manufacturer, retailer, and the entire supply chain, $i \in \{M, R, MR\}$ |
| $\alpha$ | The primary customer preference for the direct channel, $0 < \alpha < 1$ |
| $v$ | The primary customer's value coefficient for purchasing a new item |
| $u$ | The replacement customer's value coefficient for purchasing a new item |
| $s$ | The government's unit subsidy for the manufacturer's remanufacturing effort, $s > 0$ |
| $c_m, c_r$ | The unit production cost for manufacturing and remanufacturing a new item, $c_m > c_r > 0$ |
| $\Delta$ | The unit cost savings through the remanufacturing, $\Delta = c_m - c_r > 0$ |

**Assumption 3.** *The per-unit production cost for manufacturing new items $c_m$ is higher than the unit cost for remanufacturing used items $c_r$, resulting in a positive unit cost savings $\Delta$ for the manufacturer when remanufacturing recycled items, calculated as $\Delta = c_m - c_r > 0$.*

**Assumption 4.** *Customers exhibit heterogeneous valuations (v) for new items and heterogeneous costs u for returning used items, both uniformly distributed on the interval [0, 1], for primary and replacement customers, respectively.*

**Assumption 5.** *In the DCCLSC, the manufacturer and retailer share symmetrical information.*

Collectively, these assumptions encapsulate prevalent industrial realities, such as the distinct cost advantages associated with remanufacturing and the differentiation between primary and replacement customers. By elucidating the pricing and subsidy strategies employed by each stakeholder (manufacturer or retailer) within a framework of sequential decision-making, the model preserves both theoretical clarity and practical applicability. Consequently, the subsequent equilibrium analysis offers insights pertinent to typical DCCLSC operations.

## 2.3 Demand Functions and Profit Functions

In the Model M, primary customers decide whether to purchase new items from the manufacturer's direct channel or from the retailer's channel, while replacement customers trade in used items under a manufacturer-provided subsidy $b_m^M$. The utility functions are defined in Equations (1) – (3) as follow:

$$U_1^M = \alpha v - p_m^M, \tag{1}$$
$$U_2^M = v - p_r^M, \tag{2}$$
$$U_3^M = b_m^M - u, \tag{3}$$

where $\alpha \in (0, 1)$ captures the channel preference for purchasing directly from the manufacturer, $v \in [0, 1]$ denotes the valuation for purchasing a new item for primary customers, and $u \in [0, 1]$ represents the valuation of trading in used items for replacement customers.

By evaluating these utility functions, consumers make purchasing decisions based on whether their valuation exceeds the corresponding prices or subsidies. Specifically, primary customers choose the direct channel if $U_1^M \geq U_2^M$ and their valuations exceed the reservation prices; otherwise, they choose the retailer's indirect channel. Replacement customers choose to trade in used items if their valuation $u$ is lower than the offered subsidy $b_m^M$. By integrating these consumer choices over the assumed uniform distributions of $v$ and $u$, we derive the explicit demand functions for both direct and indirect channels, as presented in Equations (4)-(6). Such derivations are consistent with consumer utility maximization theory and commonly employed in closed-loop supply chain modeling (Savaskan et al., 2004). Given that and are uniformly distributed on the interval $[0, 1]$, their corresponding probability density functions are $f(v) = 1$ and $f(u)$. Using these assumptions, we obtain:

$$q_1^M = \int_{\frac{p_r^M}{\alpha}}^{\frac{p_r^M - p_m^M}{1-\alpha}} f(v)\, dv = \frac{\alpha p_r^M - p_m^M}{\alpha(1-\alpha)} \tag{4}$$

$$q_2^M = \int_{\frac{p_r^M - p_m^M}{1-\alpha}}^{1} f(v)\, dv = 1 - \frac{p_r^M - p_m^M}{1-\alpha} \tag{5}$$

$$q_3^M = \int_0^{b_m^M} f(u)\, du = b_m^M \tag{6}$$

Here, $q_1^M$ and $q_2^M$ denote the demand for new items through the direct and direct channels, respectively, while $q_3^M$ represents the quantity of used items recycled by replacement customers, incentivized by the subsidy $b_m^M$.

By integrating the demand functions into the objective functions of each member, the profit functions are consequently formulated as follows:

$$Max\ \Pi_m^M\ (p_m^M, w^M, b_m^M) = (p_m^M - c_m)q_1^M + (w^M - c_m)q_2^{MR} + (\Delta + s + p_m^M - c_m - b_m^M)q_3^M, \tag{7}$$

$$Max\ \Pi_r^M\ (p_r^M) = (p_r^M - w^M)q_2^M, \tag{8}$$

$$Max\ \Pi_s^M = \Pi_m^M + \Pi_r^M. \tag{9}$$

In the Model R, replacement customers are served through the retailer, which offers a trade-in subsidy $b_r^R$, and receives compensation via a transfer price $t^R$ from the manufacturer. The customer utility functions, as represented in Equations (10) – (12), are defined as follows:

$$U_1^R = \alpha v - p_m^R \tag{10}$$
$$U_2^R = \alpha v - p_m^R \tag{11}$$
$$U_3^R = b_r^R - \alpha u \tag{12}$$

By evaluating these utility functions, consumers determine their purchasing or trade-in decisions based on comparisons between their individual valuations and the respective channel prices or subsidies. Specifically, replacement customers opt to trade in their used items with the retailer if their valuation $u$ satisfies the condition $b_r^R \geq \alpha u$. Assuming that $v$ and $u$ follow uniform distributions within the interval $[0,1]$ with probability density functions $f(v) = 1$ and $f(u) = 1$, we integrate these utility conditions to derive the explicit demand functions, as shown in Equations (13) – (15):

$$q_1^R = \int_{\frac{p_r^R}{\alpha}}^{\frac{p_r^R - p_m^R}{1-\alpha}} f(v)\, dv = \frac{\alpha p_r^R - p_m^R}{\alpha(1-\alpha)} \tag{13}$$

$$q_2^R = \int_{\frac{p_r^R - p_m^R}{1-\alpha}}^{1} f(v)\,dv = 1 - \frac{p_r^R - p_m^R}{1-\alpha} \tag{14}$$

$$q_3^R = \int_0^{b_r^R} f(u)\,du = \frac{b_r^R}{\alpha} \tag{15}$$

Here, $q_1^R$ and $q_2^R$ represents the primary demand captured through direct and indirect channels, respectively, while $q_3^R$ denotes the recycled volume from replacement customers, incentivized by the retailer's trade-in subsidy.

These demand components are subsequently incorporated into the profit functions, as specified in Equations (16) – (18):

$$\begin{aligned} Max\ \Pi_m^R\ (p_m^R, w^R, t^R) \\ = (p_m^R - c_m)q_1^R + (w^R - c_m)q_2^R \\ + (\Delta + s + w^R - c_m - t^R)q_3^R \end{aligned} \tag{16}$$

$$Max\ \Pi_r^R\ (p_r^R, b_r^R) = (p_r^R - w^R)q_2^R + (p_r^R + t^R - w^R - b_r^R)q_3^R \tag{17}$$

$$Max\ \Pi_s^R = \Pi_m^R + \Pi_r^R \tag{18}$$

In the Model MR, both the manufacturer and the retailer concurrently provide subsidies for trade-ins, allowing replacement customers to choose either the manufacturer's subsidy $b_m^{MR}$, or the retailer's subsidy $b_r^{MR}$. Consequently, this leads to the formulation of four distinct utility functions, as represented by Equations (19) – (22):

$$U_1^{MR} = \alpha v - p_m^{MR} \tag{19}$$
$$U_2^{MR} = v - p_r^{MR} \tag{20}$$
$$U_3^{MR} = b_m^{MR} - u \tag{21}$$
$$U_4^{MR} = b_r^{MR} - \alpha u \tag{22}$$

Based on these utility comparisons, primary customers choose the direct channel when $U_1^{MR} \geq U_2^{MR}$, and replacement customers prefer the manufacturer's subsidy if $U_3^{MR} \geq U_4^{MR}$; otherwise, they opt for the retailer's subsidy. Given that valuation parameters $v$ and $u$ follow uniform distributions within the interval $[0,1]$ with probability density functions $f(v) = 1$ and $f(u) = 1$, we integrate over these intervals to derive explicit demand functions shown in Equations (23) – (26):

$$q_1^{MR} = \int_{\frac{p_r^{MR}}{\alpha}}^{\frac{p_r^{MR} - p_m^{MR}}{1-\alpha}} f(v)\, dv = \frac{\alpha p_r^{MR} - p_m^{MR}}{\alpha(1-\alpha)} \quad (23)$$

$$q_2^{MR} = \int_{\frac{p_r^{MR} - p_m^{MR}}{1-\alpha}}^{1} f(v)\,dv = 1 - \frac{p_r^{MR} - p_m^{MR}}{1-\alpha} \quad (24)$$

$$q_3^{MR} = \int_0^{\frac{b_m^{MR} - b_r^{MR}}{1-\alpha}} f(u)\, du = 1 - \frac{b_m^{MR} - b_r^{MR}}{1-\alpha} \quad (25)$$

$$q_4^{MR} = \int_{\frac{b_m^{MR} - b_r^{MR}}{1-\alpha}}^{\frac{b_r^{MR}}{\alpha}} f(u)\, du = \frac{b_r^{MR} - \alpha b_m^{MR}}{\alpha(1-\alpha)} \quad (26)$$

The profit functions are formulated as follows:

$$\begin{aligned} Max\ \Pi_m^{MR}&(p_m^{MR}, t^{MR}, b_m^{MR}, w^{MR}) \\ &= (p_m^{MR} - c_m)q_1^{MR} + (w^{MR} - c_m)q_2^{MR} \\ &\quad + (\Delta + s + p_m^{MR} - c_m - b_m^{MR})q_3^{MR} \\ &\quad + (\Delta + s + w^{MR} - c_m - t^{MR})q_4^{MR} \end{aligned} \quad (27)$$

$$\begin{aligned} Max\ \Pi_r^{MR}&(p_r^{MR}, b_r^{MR}) \\ &= (p_r^{MR} - w^{MR})q_2^{MR} + (p_r^{MR} + t^{MR} - w^{MR} - b_r^{MR})q_4^{MR} \end{aligned} \quad (28)$$

$$Max\ \Pi_s^{MR} = \Pi_m^{MR} + \Pi_r^{MR} \quad (29)$$

In this framework, $q_1^{MR}$ represents the primary demand directed through the manufacturer, while $q_2^{MR}$ denotes the primary demand routed through the retailer. Furthermore, $q_3^{MR}$ signifies the segment of replacement demand that favors the manufacturer's subsidy, and $q_4^{MR}$ indicates the segment that prefers the retailer's subsidy.

By solving these maximization problems, we can determine the equilibrium prices, and recycling subsidies that align the objectives of manufacturers and retailers within each dual-channel structure. This ensures optimal performance in DCCLSCs, balancing economic and environmental goals.

## 3  Model Solution and Analysis

In this section, we derive the equilibrium solutions for Model M, Model R, and Model MR. The manufacturer, acting as the Stackelberg leader, determines the wholesale price, direct selling price, and recycling subsidy, while the retailer, as the follower, determines the retail selling price and recycling subsidy. For each model, we solve the profit-maximization problem by deriving the first-order conditions (FOCs) for both the manufacturer and the retailer, and then iteratively solving the system of equations. The Hessian matrix analysis is employed to confirm the uniqueness and stability of the equilibrium solutions.

### 3.1 Equilibrium Results for Model M

In Model M, the manufacturer is responsible for collecting used items, which significantly influences pricing decisions in both the forward and reverse supply chains. The equilibrium is determined using a backward induction approach. First, we derive the FOC for the selling price in the indirect channel:

$$\frac{\partial \Pi_r^M}{\partial p_r^M} = q_2^M + (p_r^M - w^M)\frac{\partial \Pi q_2^M}{\partial p_r^M} = 0 \tag{30}$$

Substituting the expression for $q_2^M$ and solving for $p_r^M$ in terms of $w^N$ and $p_m^M$:

$$p_r^{M*} = \frac{w^{N*} + (1-\alpha)p_m^{M*}}{2-\alpha} \tag{31}$$

This yields the retailer's optimal selling price in response to the manufacturer's decision variables. Then, we derive the FOCs for the manufacturer's decision variables:

$$\frac{\partial \Pi_m^M}{\partial p_m^M} = 0 \tag{32}$$

$$\frac{\partial \Pi_m^M}{\partial w^M} = 0 \tag{33}$$

$$\frac{\partial \Pi_m^M}{\partial b_m^M} = 0 \tag{34}$$

Next, substituting the retailer's optimal selling price $p_r^{M*}$ into the manufacturer's FOCs and solving for $p_m^M$, $w^M$, and $b_m^M$, we obtain the following equilibrium solutions:

$$p_m^{M*} = -\frac{2\alpha + 2c_m + \Delta\alpha - \alpha c_m + \alpha s}{\alpha - 4} \tag{35}$$

$$w^{M*} = -\frac{4c_m - \alpha + 2\Delta\alpha - 2\alpha c_m + 2\alpha s + \alpha^2 + 4}{2(\alpha - 4)} \tag{36}$$

$$b_m^{M*} = -\frac{2\Delta + \alpha - c_m + 2s}{\alpha - 4} \tag{37}$$

Finally, substituting these equilibrium solutions into the retailer's pricing reaction function, we solve for the retailer's optimal selling price:

$$p_r^{M*} = -\frac{8c_m - 7\alpha + 4\Delta\alpha - 4\alpha c_m + 4\alpha s + 3\alpha^2 + 12}{4(\alpha - 4)} \tag{38}$$

**Theorem 1.** *Within the framework of the manufacturer-led recycling framework, a unique equilibrium solution exists for the unit selling price of new items in the manufacturer's direct channel $p_m^{M^*}$, the unit wholesale price $w^{M^*}$, the unit recycling subsidy for used items $b_m^{M^*}$, and the unit selling price in the retailer's indirect channel $p_r^{M^*}$.*

**Proposition 1.** *In the Model M, the equilibrium unit selling price in the retailer's indirect channel is lower than that in the manufacturer's direct channel, i.e., $p_m^{M^*} < p_r^{M^*}$.*

*Proof.* From the optimal solutions in the Model M, we compute the difference between the direct and indirect selling prices :

$$p_m^{M^*} - p_r^{M^*} = \frac{-4(2\alpha + 2c_m + \Delta\alpha - \alpha c_m + \alpha s) + (8c_m - 7\alpha + 4\Delta\alpha - 4\alpha c_m + 4\alpha s + 3\alpha^2 + 12)}{4(\alpha - 4)} \quad (39)$$

which simplifies to:

$$p_m^{M^*} - p_r^{M^*} = \frac{3\alpha^2 - 15\alpha + 12}{4(\alpha - 4)} \quad (40)$$

The expression in the equation (40) exhibits a negative denominator due to $0 < \alpha < 1$, which implies $\alpha - 4 < 0$. Upon simplification, the numerator is represented by the quadratic function $3\alpha^2 - 15\alpha + 12$. This function is positive within the interval $0 < \alpha < 1$ because its roots are located at $\alpha = 1$ and $\alpha = 4$, and the probola opens upwards. Consequently, a positive numerator divided by a negative denominator yields a negative value, indicating that the entire expression is less than zero for $0 < \alpha < 1$. Therefore, it follows that $p_m^{M^*} < p_r^{M^*}$.

Proposition 1 elucidates that within a manufacturer-led recycling framework, the manufacturer strategically sets a lower direct selling price than the retailer. This outcome stems from the Stackelberg game framework, where the manufacturer, as the leader, reduces its direct price to capture demand and enhance the appeal of its trade-in program. In contrast, the retailer, which does not directly benefit from remanufacturing, is compelled to set a higher price to maintain profitability. These findings highlight the critical role of customer heterogeneity: primary customers, who purchase new items without trading in used items, face higher retail prices, while replacement customers, incentivized by the manufacturer's subsidy, benefit from reduced direct selling prices. This insight underscores the importance of pricing coordination in DCCLSCs, where manufacturers leverage direct sales pricing to drive recycling participation and promote overall supply chain sustainability.

**Proposition 2.** *In the Model M, the equilibrium decision and their dependence on the primary customer preference for the direct channel $\alpha$ are characterized as follow:*

i. *The optimal wholesale price $w^{M^*}$ increases in $\alpha$ if and only if $c_m < 2\Delta + 2s + 1$. Otherwise.*

ii. *The unit subsidy specified by the manufacturer $b_m^{M^*}$, and the selling price in the direct channel $p_m^{M^*}$ increase in $\alpha$ if and only if $c_m < 2\Delta + 2s + 4$. Otherwise.*

iii. The selling price in the indirect channel $p_r^{M^*}$ increases in $\alpha$ if and only if $c_m < \frac{4(\Delta+s)-1}{2}$.

*Otherwise.*

*Proof.* Note that each decision variable can be written in the form:

$$X(\alpha) = -\frac{F(\alpha)}{k(\alpha-4)}, \tag{41}$$

where $k > 0$ is a constant (possibly 1, 2, or 4) and $F(\alpha)$ is a polynomial in $\alpha$ that also includes parameters $\{c_m, c_r, s, \Delta\}$. Since $0 < \alpha < 1$ implies $\alpha - 4 < 0$, the sign of the derivate:

$$\frac{\partial X}{\partial \alpha} = -\frac{F'(\alpha)k(\alpha-4) - kF(\alpha)}{k^2(\alpha-4)^2} = -\frac{F'(\alpha)(\alpha-4) - F(\alpha)}{k(\alpha-4)^2}, \tag{42}$$

depends only on:

$$F'(\alpha)(\alpha-4) - F(\alpha). \tag{43}$$

Under normal cost-parameter ranges, $F'(\alpha)(\alpha-4) - F(\alpha)$ does not switch sign for $\alpha \in (0,1)$, so $X(\alpha)$ is strictly monotonic. Consequently, comparing the endpoint value $X(0)$ and $X(1)$ suffices to determine whether $X$ increases or decreases in $\alpha$:

i. If $X(1) > X(0)$, then $X(\alpha)$ increases on $(0,1)$.
ii. If $X(1) < X(0)$, then $X(\alpha)$ decreases on $(0,1)$.

Thus, each threshold condition arises from:

$$w^{M^*}(0) = \frac{c_m+1}{2}, \; w^{M^*}(1) = \frac{c_m+\Delta+s+2}{3} \Rightarrow w^{M^*}(1) > w^{M^*}(0)$$

$$\Rightarrow c_m < 2\Delta + 2s + 1,$$

$$b_m^{M^*}(0) = \frac{2\Delta-c_m+2s}{4}, \; b_m^{M^*}(1) = \frac{2\Delta-c_m+2s-1}{3} \Rightarrow b_m^{M^*}(1) > b_m^{M^*}(0)$$

$$\Rightarrow c_m < 2\Delta + 2s + 4, \tag{44}$$

$$p_m^{M^*}(0) = \frac{c_m}{2}, \; p_m^{M^*}(1) = \frac{2+c_m+\Delta+s}{3} \Rightarrow p_m^{M^*}(1) > p_m^{M^*}(0)$$

$$\Rightarrow c_m < 2\Delta + 2s + 4,$$

$$p_r^{M^*}(0) = \frac{2c_m+3}{4}, \; p_r^{M^*}(1) = \frac{c_m+\Delta+s+2}{3} \Rightarrow p_r^{M^*}(1) > p_r^{M^*}(0),$$

$$\Rightarrow c_m < \frac{4(\Delta+s)-1}{2}.$$

Proposition 2 further demonstrates that within a manufacturer-led recycling framework, an increasing propensity among primary customers to purchase through the direct channel can have dual effects on pricing dynamics. Specifically, it can either elevate or depress the wholesale price, the per-unit subsidy (buyback) offered to replacement customers trading in used items, and the final selling prices across both channels. When the cost of producing new goods is sufficiently low relative to the manufacturer's potential gains from material reuse and governmental support, heightened direct-channel demand leads to higher subsidies for replacement customers and enhanced pricing power. Conversely, if manufacturing and remanufacturing costs become excessively burdensome, the same increase in direct-channel demand necessitates downward price adjustments and reduces the buyback level for used products. Thus, operational efficiency, effective remanufacturing practices, and strategic utilization of policy incentives are essential for converting increased direct-channel appeal among primary customers into profitable outcomes while maintaining attractive trade-in subsidies for replacement customers.

This section demonstrates that manufacturers operating within a manufacturer-led recycling framework can effectively implement differentiated pricing strategies to cater to both primary and replacement customers. By offering competitive direct selling prices and trade-in subsidies, manufacturers can simultaneously stimulate demand from primary customers and encourage the return of used items from replacement customers. However, the success of these strategies is contingent on the manufacturer's ability to maintain efficient production processes and leverage remanufacturing opportunities. For instance, in the electronics industry, companies that design modular products can reduce disassembly labor and take advantage of government incentives, thereby supporting sustainable practices. This approach not only facilitates higher pricing and increased trade-in subsidies but also maintains healthy profit margins. Managers should focus on optimizing production processes and seeking policy support to overcome potential challenges and maximize the benefits of their pricing and recycling strategies. This alignment ensures both economic viability and environmental sustainability, driving progress toward a more circular and sustainable supply chain model.

## 3.2 Equilibrium Results for Model R

Presented below are the distinct equilibrium solutions for the Model R, where the retailer is responsible for collecting used items, as outlined in Equations (45)-(49). In contrast to Model M, the manufacturer is required to offer the retailer a per-unit transfer price $t^R$ to compensate for the retailer's recycling efforts. These results are also obtained by applying backward induction to the profit maximization problems of both the retailer and the manufacturer, followed by the simultaneous resolution of each decision variable. The resulting closed-form expressions describe the equilibrium in terms of the primary customer's preference for the direct channel $\alpha$, the unit production cost for remanufacturing a new item $c_m$, the unit cost savings through the remanufacturing $\Delta$, and the government's unit subsidy for the manufacturer's remanufacturing effort $s$.

$$p_m^{R^*} = \frac{2\Delta\alpha - 2c_m - \alpha + 8\alpha c_m + 2\alpha s + 9\alpha^2}{18\alpha - 4} \tag{45}$$

$$w^{R^*} = \frac{5\alpha - c_m + \Delta\alpha + 4\alpha c_m + \alpha s - 1}{9\alpha - 2} \tag{46}$$

$$b_r^{R^*} = \frac{\alpha(2\Delta - c_m + 2s + 1)}{9\alpha - 2} \tag{47}$$

$$p_r^{R^*} = \frac{4\Delta + 29\alpha - 6c_m + 4s + 18\alpha c_m - 9\alpha^2 - 4}{36\alpha - 8} \tag{48}$$

$$t^{R^*} = -\frac{4\Delta + \alpha - 2c_m + 4s - 20\Delta\alpha + 10\alpha c_m - 20\alpha s - 9\alpha^2}{4(9\alpha - 2)} \tag{49}$$

**Theorem 2.** *Within the retailer-led recycling framework, there exists a unique equilibrium solution for the unit selling price of new items in the manufacturer's direct channel $p_m^{R^*}$, the unit wholesale price $w^{R^*}$, the unit recycling subsidy for used items $b_r^{R^*}$, the unit selling price in the retailer's indirect channel $p_r^{R^*}$, and the unit transfer price for the retailer specified by the manufacturer for collecting used items $t^{R^*}$.*

**Proposition 3.** *In the Model R, the relationship between the equilibrium unit selling prices for the indirect channel $p_r^{R^*}$ and the direct channel $p_m^{R^*}$ depends on the primary customer preference for the direct channel $\alpha$:*

i.     *If $\frac{2}{9} < \alpha < 1$, then $p_m^{R^*} < p_r^{R^*}$.*

ii.     *If $0 < \alpha < \frac{2}{9}$, then $p_m^{R^*} > p_r^{R^*}$.*

*Proof.* To prove this, we compute the difference between the two equilibrium prices:

$$\begin{aligned} D &= p_m^{R^*} - p_r^{R^*} \\ &= \frac{4\Delta(\alpha - 1) + 2c_m - 31\alpha - 2\alpha c_m + 4\alpha s - 4s + 27\alpha^2 + 4}{36\alpha - 8} \end{aligned} \tag{50}$$

To determine the sign of $D$, we focus on the numerator:

$$\begin{aligned} N &= 4\Delta(\alpha - 1) + 2c_m - 31\alpha - 2\alpha c_m + 4\alpha s - 4s + 27\alpha^2 + 4 \\ &= 27\alpha^2 + (4\Delta - 31 - 2c_m + 4s)\alpha \\ &\quad + (2c_m - 4\Delta - 4s + 4) \end{aligned} \tag{51}$$

Setting $D = 0$, we solve for $\alpha$ by considering the quadratic equation:

$$27\alpha^2 + B\alpha + C = 0, \tag{52}$$

where

$$B = 4\Delta - 31 - 2c_m + 4s, \quad C = 2c_m - 4\Delta - 4s + 4 \tag{53}$$

Applying the quadratic formula:

$$\alpha = \frac{-B \pm \sqrt{B^2 - 108C}}{54} \tag{54}$$

Under specific parameter conditions, one of the roots simplifies to $\frac{2}{9}$, which serves as a threshold value. The sign of $D$ is determined as follow:

$$\text{If } \frac{2}{9} < \alpha < 1, \text{ then } D < 0 \Rightarrow p_m^{R^*} < p_r^{R^*},$$

$$\text{If } 0 < \alpha < \frac{2}{9}, \text{ then } D > 0 \Rightarrow p_m^{R^*} > p_r^{R^*}. \tag{55}$$

Proposition 3 demonstrates that the preferences of primary customers for direct purchasing significantly influence the relative pricing strategies employed by manufacturers in their direct sales channels and retailers in their indirect sales channels, particularly within a retailer-led recycling framework where the retailer governs the trade-in policy. In response to an increased inclination among primary customers to purchase directly from the manufacturer, the manufacturer strategically reduces its direct selling price to attract a larger customer base. Conversely, the retailer experiences diminished demand within its channel and subsequently raises its selling price to maintain profitability. This dynamic allows the manufacturer to enhance the attractiveness of direct purchases, thereby compensating for its limited control over trade-in incentives. Conversely, when the retailer's channel is favored by its primary customers, the retailer may further solidify its competitive advantage by decreasing its prices. In this scenario, the manufacturer, faced with a reduction in demand for its direct channel, may respond by increasing its selling price to optimize revenue from a shrinking customer base.

**Proposition 4.** *In the Model R, the equilibrium decisions and their dependence on the primary customer preference for the direct channel $\alpha$ are characterized as follows:*

i. *The optimal wholesale price $w^{R^*}$ increases in $\alpha$ if and only if $c_m < \frac{5+\Delta+s}{4}$. Otherwise.*

ii. *The unit subsidy specified by the retailer $b_r^{R^*}$ increases in $\alpha$ if and only if $c_m < 2\Delta + 2s + 1$. Otherwise.*

iii. *The selling price in the direct channel $p_m^{R^*}$ increase in $\alpha$ if and only if $c_m < \frac{2\Delta+2s+8}{3}$. Otherwise.*

iv. *The selling price in the indirect channel $p_r^{R^*}$ increases in $\alpha$ if and only if $c_m < \frac{4(\Delta+s)-1}{2}$. Otherwise.*

v. *The unit transfer price for the retailer specified by the manufacturer for collecting used items $t^{R^*}$ increases in $\alpha$ if and only if $c_m > \frac{4(\Delta+s)+9}{2}$. Otherwise.*

*Proof.* Each equilibrium variable $X(\alpha)$ can be expressed in the general form:

$$X(\alpha) = \frac{N(\alpha)}{D(\alpha)}, \tag{56}$$

where $N(\alpha)$ is a polynomial in $\alpha$, incorporating parameters $c_m$, $\Delta$, $s$. $D(\alpha)$ is a function of $\alpha$, such as $18\alpha - 4$, $9\alpha - 2$, etc.

The derivative of $X(\alpha)$ with respect to $\alpha$ is given by:

$$\frac{\partial X}{\partial \alpha} = \frac{N'(\alpha)D(\alpha) - N(\alpha)D'(\alpha)}{D(\alpha)^2}. \tag{57}$$

Since $D(\alpha)^2 > 0$ for $0 < \alpha < 1$, the monotonicity of $X(\alpha)$ is determined by the sign of the numerator: $N'(\alpha)D(\alpha) - N(\alpha)D'(\alpha)$.

Thus, each threshold condition arises from:

$$w^{R^*}(0) = \frac{c_m - 1}{2}, \ w^{R^*}(1) = \frac{5 + c_m + \Delta + s}{4} \Rightarrow w^{R^*}(1) > w^{R^*}(0)$$

$$\Rightarrow c_m < \frac{5 + \Delta + s}{4},$$

$$b_r^{R^*}(0) = 0, \ b_R^{R^*}(1) = \frac{2\Delta + 2s + 1 - c_m}{7} \Rightarrow b_r^{R^*}(1) > b_r^{R^*}(0)$$

$$\Rightarrow c_m < \frac{2\Delta + 2s + 8}{3},$$

$$p_m^{R^*}(0) = \frac{c_m}{2}, \ p_m^{R^*}(1) = \frac{2\Delta + 2s + 8 + c_m}{14} \Rightarrow p_m^{R^*}(1) > p_m^{R^*}(0) \tag{58}$$

$$\Rightarrow c_m < 2\Delta + 2s + 4,$$

$$p_r^{R^*}(0) = \frac{2c_m + 3}{4}, \ p_r^{R^*}(1) = \frac{c_m + \Delta + s + 2}{3} \Rightarrow p_r^{R^*}(1) > p_r^{R^*}(0),$$

$$\Rightarrow c_m < \frac{4(\Delta + s) - 1}{2},$$

$$t^{R^*}(0) = \frac{4\Delta + 4s - 2c_m + 9}{4}, \ t^{R^*}(1) = \frac{20\Delta\alpha - 10\alpha c_m + 20\alpha s + 9\alpha^2}{4} \Rightarrow t^{R^*}(1) > t^{R^*}(0),$$

$$\Rightarrow c_m > \frac{4(\Delta + s) + 9}{2}.$$

Proposition 4 demonstrates that within a retailer-led recycling framework, the primary customer preference for direct purchasing not only influences pricing strategies but also dictates the structuring of trade-in incentives by both the manufacturer and the retailer. A critical insight is that the retailer, as the leader in recycling initiatives, exercises more direct control over trade-in subsidies, whereas the manufacturer's influence is primarily exerted through wholesale pricing and transfer payments for collected used items. When primary customers exhibit a strong preference for direct purchasing, the manufacturer benefits from heightened demand, enabling it to increase wholesale prices while maintaining competitive direct selling prices. Consequently, the retailer is compelled to enhance its trade-in subsidy to remain appealing to replacement customers. However, if the manufacturer's production costs are excessively high, its capacity to elevate wholesale and direct selling prices becomes limited, necessitating a greater reliance on adjusting transfer payments to the retailer. Conversely, when primary customers favor the retailer's channel, the retailer acquires increased control over pricing, and retailer-led trade-in incentives become a more significant driver of customer decision-making.

This section emphasizes the importance of trade-in program management in shaping pricing strategies within DCCLSCs. Firms managing recycling programs enjoy more pricing flexibility, while others must adapt to stay competitive. Manufacturers can counter retailer trade-in incentives with exclusive direct sales promotions, and retailers can adjust trade-in subsidies to optimize demand across channels. Cost efficiency is crucial for pricing adaptability; low production costs allow aggressive pricing for market share, while rising costs require careful transfer price adjustments to maintain profits. Proactive cost control and dynamic pricing are vital for maximizing profitability and promoting sustainable product recovery.

### 3.3 Equilibrium Results for Model MR

Below are the distinct equilibrium solutions for the Model MR, in which the manufacturer and the retailer collaborate on collecting used items. These results are derived by applying backward induction to the profit maximization problems of both the retailer and the manufacturer, solving the decision variables simultaneously. The resulting expressions characterize the equilibrium as follows:

$$p_m^{MR*} = -\frac{X_1 - 2\alpha + 12\alpha^2 - 6\alpha^3}{3\alpha^2 - 17\alpha + 2\alpha^3 + 4} \tag{59}$$

$$w^{MR*} = -\frac{17\alpha + 2X_1 + 4\alpha^2 - 11\alpha^3 + 2\alpha^4 - 4}{2(3\alpha^2 - 17\alpha + 2\alpha^3 + 4)} \tag{60}$$

$$b_m^{MR*} = \frac{2X_2 + X_3 - 19\alpha + 8\Delta\alpha - 18\alpha c_m + 18\alpha s - 4\alpha^2 + 11\alpha^3 + 4}{2(3\alpha^2 - 17\alpha + 2\alpha^3 + 4)} \tag{61}$$

$$b_r^{MR*} = \frac{\alpha(5\alpha + X_2 - X_3 + 3\alpha^2 - 4\alpha^3)}{3\alpha^2 - 17\alpha + 2\alpha^3 + 4} \tag{62}$$

$$p_r^{MR*} = \frac{4\Delta + 23\alpha + 4s + X_1 + \Delta\alpha - 5\alpha^2 - 9\alpha^3 + 3\alpha^4 - 4}{2(3\alpha^2 - 17\alpha + 2\alpha^3 + 4)} \tag{63}$$

$$t^{MR*} = \frac{(\alpha + 1)(\alpha + X_2 + X_3 - 6\alpha^2 + 3\alpha^3)}{3\alpha^2 - 17\alpha + 2\alpha^3 + 4} \tag{64}$$

From the equation (58) to (63), we define $X_1$, $X_2$, and $X_3$ as follows:

$$X_1 = 3\Delta\alpha - 2c_m + 7\alpha c_m + 3\alpha s - 5\Delta\alpha^2 + 2\Delta\alpha^3 + \alpha^2 c_m - 2\alpha^3 c_m \\ - 5\alpha^2 s + 2\alpha^3 \quad (65)$$

$$X_2 = 2\Delta - c_m + 2s - 10\Delta\alpha + 5\alpha c_m \quad (66)$$

$$X_3 = -10\alpha s + 4\Delta\alpha^2 - 2\alpha^2 c_m + 4\alpha^2 s \quad (67)$$

**Theorem 3.** *Within the Model MR, there exists a unique equilibrium solution for the unit selling price of new items in the manufacturer's direct channel $p_m^{MR^*}$, the unit wholesale price $w^{MR^*}$, the unit recycling subsidy specified by the manufacturer for used items $b_m^{MR^*}$, the unit recycling subsidy specified by the retailer for used items $b_r^{MR^*}$, the unit selling price in the retailer's indirect channel $p_r^{MR^*}$, and the unit transfer price for the retailer specified by the manufacturer for collecting used items $t^{MR^*}$.*

**Proposition 5.** *In the Model MR, the relationship between the equilibrium unit selling prices for the indirect channel $p_r^{MR^*}$ and the direct channel $p_m^{MR^*}$ depend on the primary customer preference for the direct channel $\alpha$. There exists a threshold $\alpha^*$, given by:*

$$\alpha^* = \frac{6c_m - 4\Delta - 4s + 4}{3c_m - 3\Delta - 3s + 21} \quad (68)$$

such that:

i.   If $\alpha > \alpha^*$, then $p_m^{MR^*} < p_r^{MR^*}$.
ii.  If $\alpha < \alpha^*$, then $p_m^{MR^*} > p_r^{MR^*}$.
iii. If $\alpha = \alpha^*$, then $p_m^{MR^*} = p_r^{MR^*}$.

*Proof.* Setting $p_m^{MR^*} = p_r^{MR^*}$ and clearing denominators, we obtain:

$$-2(X_1 - 2\alpha + 12\alpha^2 - 6\alpha^3) \\ = 4\Delta + 23\alpha + 4s + X_1 + \Delta\alpha - 5\alpha^2 - 9\alpha^3 + 3\alpha^4 \quad (69) \\ - 4.$$

Solving for $\alpha^*$, we can obtain the expression in Equation (67). In addition, if $\alpha > \alpha^*$, the numerator of $p_m^{MR^*} - p_r^{MR^*}$ becomes negative, implying $p_m^{MR^*} < p_r^{MR^*}$. Otherwise.

Proposition 5 posits that within a collaborative recycling framework, where both firms engage in the collection of used items and implement trade-in programs, the primary consumer preference for direct purchasing influences the price competition between channels. In contrast to recycling models led by either retailers or manufacturers, both entities now vie for replacement customers, rendering pricing decisions more mutually dependent. When consumers exhibit a strong preference for purchasing directly from the manufacturer, the manufacturer reduces its direct selling price to attract demand, while the retailer maintains a higher price due to diminished competitiveness in the indirect channel. Conversely, when customers favor the retailer's channel, the retailer decreases its price to enhance its competitive position, prompting the manufacturer to increase its direct selling price in response to the diminished demand.

**Proposition 6.** *In the Model MR, the unit recycling subsidy specified by the manufacturer for used items $b_m^{MR^*}$ and the unit recycling subsidy specified by the retailer for used items $b_r^{MR^*}$ depend on the primary customer preference for the direct channel $\alpha$. There exists a threshold $\dot{\alpha}$, given by:*

$$\dot{\alpha} = \frac{5c_m - 10\Delta - 10s + 8}{10\Delta - 5c_m + 10s + 2} \tag{70}$$

such that:

i. If $\alpha > \dot{\alpha}$, then $b_m^{MR^*} > b_r^{MR^*}$.
ii. If $\alpha < \dot{\alpha}$, then $b_m^{MR^*} < b_r^{MR^*}$.
iii. If $\alpha = \dot{\alpha}$, then $b_m^{MR^*} = b_r^{MR^*}$.

*Proof.* To compare the subsidies, we define the difference:

$$D = b_m^{MR^*} - b_r^{MR^*} \tag{71}$$

Substituting their equilibrium expressions:

$$D = \frac{2X_2 + X_3 - 19\alpha + 8\Delta\alpha - 18\alpha c_m + 18\alpha s - 4\alpha^2 + 11\alpha^3 + 4}{2(3\alpha^2 - 17\alpha + 2\alpha^3 + 4)}$$

$$- \frac{\alpha(5\alpha + X_2 - X_3 + 3\alpha^2 - 4\alpha^3)}{3\alpha^2 - 17\alpha + 2\alpha^3 + 4} \tag{72}$$

Solving for $\dot{\alpha}$, we obtain the result in Equation (69). If $\alpha > \dot{\alpha}$, the numerator of $b_m^{MR^*} - b_r^{MR^*}$ becomes positive, implying $b_m^{MR^*} > b_r^{MR^*}$. Otherwise.

Proposition 6 indicates that both manufacturers and retailers offer trade-in subsidies, with the higher subsidy depending on primary customer preference for direct purchasing. If primary customers prefer purchasing directly from the manufacturer, the manufacturer offers a higher subsidy to attract more replacement customers, while the retailer reduces its subsidy due to lower demand. Conversely, if primary customers prefer the retailer's channel, the retailer offers a higher subsidy to attract more used items, and the manufacturer lowers its subsidy to manage costs. This underscores the importance of strategic coordination over competitive subsidy increases. By balancing trade-in incentives and optimizing subsidies based on cost structures, manufacturers and retailers can enhance profitability.

**Proposition 7.** *In the Model MR, the equilibrium decisions and their dependence on the primary customer preference for the direct channel $\alpha$ are characterized as follows:*

i. *The optimal wholesale price $w^{MR^*}$ increases in $\alpha$ if and only if $c_m < \frac{10\Delta + 6s + 8}{9}$. Otherwise.*

ii. *The unit subsidy specified by the manufacturer $b_m^{MR^*}$ increases in $\alpha$ if and only if $c_m < 2\Delta + 2s + 1$. Otherwise. While the unit subsidy specified by the retailer $b_r^{MR^*}$ increases in $\alpha$ if and only if $c_m < \frac{8\Delta + 8s + 5}{5}$. Otherwise.*

iii. *The selling price in the direct channel $p_m^{MR^*}$ increase in $\alpha$ if and only if $c_m < \frac{6\Delta + 6s + 4}{7}$. Otherwise. While the selling price in the indirect channel $p_r^{MR^*}$ increases in $\alpha$ if and only if $c_m < \frac{6\Delta + 6s + 4}{5}$. Otherwise.*

iv. *The unit transfer price for the retailer specified by the manufacturer for collecting used items $t^{MR^*}$ increases in $\alpha$ if and only if $c_m > \frac{6\Delta + 6s + 9}{4}$. Otherwise.*

*Proof.* For each decision variable $X(\alpha)$, we analyze its monotonicity by computing its derivate:

$$\frac{dX}{d\alpha} = \frac{N'(\alpha)D(\alpha) - N(\alpha)D'(\alpha)}{D(\alpha)^2} \tag{73}$$

Since $D(\alpha)^2$ is always positive, the sign of $\frac{dX}{d\alpha}$ depends on:

$$N'(\alpha)D(\alpha) - N(\alpha)D'(\alpha) \tag{74}$$

Thus,

$$\frac{dp_m^{MR^*}}{d\alpha} > 0 \Rightarrow c_m < \frac{6\Delta + 6s + 4}{7}, \quad \frac{dw^{MR^*}}{d\alpha} > 0 \Rightarrow c_m < \frac{10\Delta + 6s + 8}{9},$$

$$\frac{db_m^{MR^*}}{d\alpha} > 0 \Rightarrow c_m < 2\Delta + 2s + 1, \quad \frac{db_r^{MR^*}}{d\alpha} > 0 \Rightarrow c_m < \frac{8\Delta + 8s + 5}{5}, \tag{75}$$

$$\frac{dp_r^{MR^*}}{d\alpha} > 0 \Rightarrow c_m < \frac{6\Delta + 6s + 4}{5}, \quad \frac{dt^{MR^*}}{d\alpha} < 0 \Rightarrow c_m > \frac{6\Delta + 6s + 9}{4}.$$

Proposition 7 elucidates that within Model MR, the wholesale price, the manufacturer's unit subsidy, and the selling price in the direct channel all escalate with an intensified primary customer preference for the direct channel, provided that the manufacturing cost remains below specific thresholds. Concurrently, the retailer's trade-in subsidy and the selling price in the indirect channel also ascend with a stronger preference for the direct channel, contingent upon their respective cost conditions. However, when the manufacturing cost exceeds a critical threshold, the unit transfer price for used-item collection increases, thereby altering the incentive structure. In situations where primary customers demonstrate a stronger preference for the direct channel and the manufacturer's costs are sufficiently low, the manufacturer is positioned to elevate its wholesale price, direct channel selling price, and trade-in subsidy, prompting the retailer to adjust its strategy accordingly. Conversely, if costs rise or customer preference shifts towards the retailer, the retailer is afforded the opportunity to enhance its own subsidy and selling price. Understanding these cost-dependent tipping points enables both firms to effectively coordinate their decisions, ensuring profitability and maintaining efficient used-item recovery within Model MR.

These results provide key insights for managing collaborative recycling. When manufacturers and retailers jointly recover used products, the pricing gap between direct and indirect channels hinges on customer preference for the direct channel. If customers strongly prefer the direct channel, manufacturers lower direct prices and increase recycling subsidies, leading retailers to raise prices and reduce subsidies.

# 4 Numerical Analysis

To validate the results of the equilibrium analysis, a sensitivity analysis was conducted on key decision variables across various parameter settings, as detailed in **TABLE 4**. For Model M, the parameter $\alpha$ (representing primary customers' preference for the direct channel) was set at 0.70 and 0.80, while maintaining $c_m = 1.2$, $c_r = 1.0$, and $s$ (a governmental unit subsidy for the manufacturer's remanufacturing effort) at either 0.10 or 0.15. In the case of Model R, $\alpha$ was selected as 0.65 and 0.75, with $c_m = 1.5$, $c_r = 0.7$ and $s$ at 0.20 or 0.25. Finally, for Model MR, $\alpha$ was set at 0.60 and 0.70, with $c_m = 1.0$, $c_r = 0.5$. and $s$ at 0.20 or 0.30.

**TABLE 4.** Numerical Analysis

| Model | $\alpha$ | $c_m$ | $c_r$ | $\Delta$ | $s$ | $p_m^{i*}$ | $p_r^{i*}$ | $w^{i*}$ | $b_m^{i*}$ | $b_r^{i*}$ | $t^{i*}$ |
|---|---|---|---|---|---|---|---|---|---|---|---|
| M | 0.7 | 1.2 | 1 | 0.2 | 0.1 | 1.3 | 1.6 | 1.15 | 0.4 | — | — |
| M | 0.8 | 1.2 | 1 | 0.2 | 0.15 | 1.4 | 1.7 | 1.25 | 0.45 | — | — |
| R | 0.65 | 1.5 | 0.7 | 0.8 | 0.2 | 1.1 | 1.5 | 1 | — | 0.5 | 0.4 |
| R | 0.75 | 1.5 | 0.7 | 0.8 | 0.25 | 1.3 | 1.7 | 1.1 | — | 0.55 | 0.45 |
| MR | 0.6 | 1 | 0.5 | 0.5 | 0.2 | 1 | 1.4 | 1.2 | 0.35 | 0.25 | 0.5 |
| MR | 0.7 | 1 | 0.5 | 0.5 | 0.3 | 1.2 | 1.6 | 1.3 | 0.4 | 0.3 | 0.55 |

Analyzing the impact of varying $\alpha$ values (or adjustments in $\{c_m, c_r, s\}$ on decision-making reveals an obvious pattern: an increase in $\alpha$ typically results in the manufacturer adopting a more aggressive pricing strategy, while the retailer's selling price may remain relatively stable if $\Delta$ or the cost structure is not excessively large. This observation is consistent with Proposition 1 in Model M, which posits that $p_m^{M*} < p_r^{M*}$, indicating that the unit selling price in the indirect channel is lower than in the direct channel. The rationale behind this is that, despite receiving a government subsidy $s$ to promote remanufacturing, the manufacturer maintains a higher price in the direct channel while ensuring sufficient margin for the retailer, thereby minimizing potential conflicts.

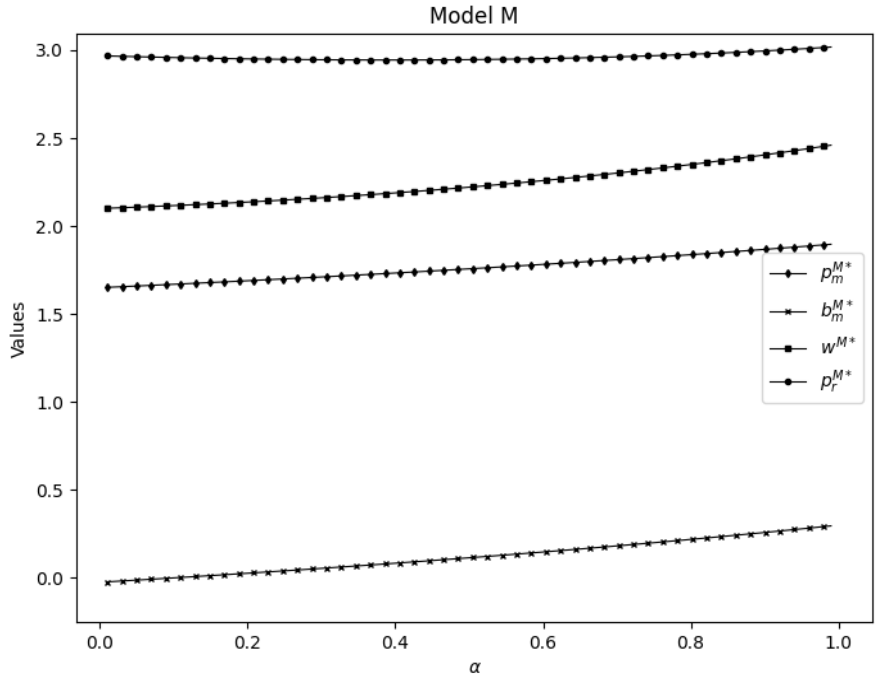

**FIGURE 3**. Decision Variables for Model M

**FIGURE 3** illustrates the decision variables for Model M across a spectrum of $\alpha$ values ranging from 0.01 to 0.99. The simulation parameters are set as follows: $s = 1.5$, $c_m = 6$, $c_r = 4$, and $\Delta = 2$. The decision variables depicted include $p_m^{M*}$, $p_r^{M*}$, $w^{M*}$, and $b_m^{M*}$. The findings reveal a monotonic increase in all variables with rising $\alpha$, highlighting a strong correlation between the primary customer's preference for the direct channel and price levels. Notably, the relationship $p_r^{M*} > p_m^{M*}$ is maintained, indicating a clear pricing hierarchy. The subsidy $b_m^{M*}$ exhibits an upward trend, suggesting enhanced manufacturer involvement in incentivizing used products as the primary customer's preference for the direct channel increases. Based on these findings, a manufacturer-led recycling framework empowers manufacturers to strategically modify pricing and trade-in incentives in response to primary customer loyalty towards direct sales channels. For instance, an electronics manufacturer experiencing an increase in online sales may choose to raise prices within direct channels while simultaneously enhancing subsidies for replacement customers. This strategy ensures a consistent supply of used items for remanufacturing. By adopting such an integrated approach, firms not only maintain profitability but also enhance circular operations through the reliable procurement of trade-in items essential for the remanufacturing process. Moreover, establishing the retail price above the direct-channel price aligns with theoretical insights, which suggest that a pronounced customer preference for direct channels allows manufacturers to influence the overall pricing structure while safeguarding the retailer's profit margins.

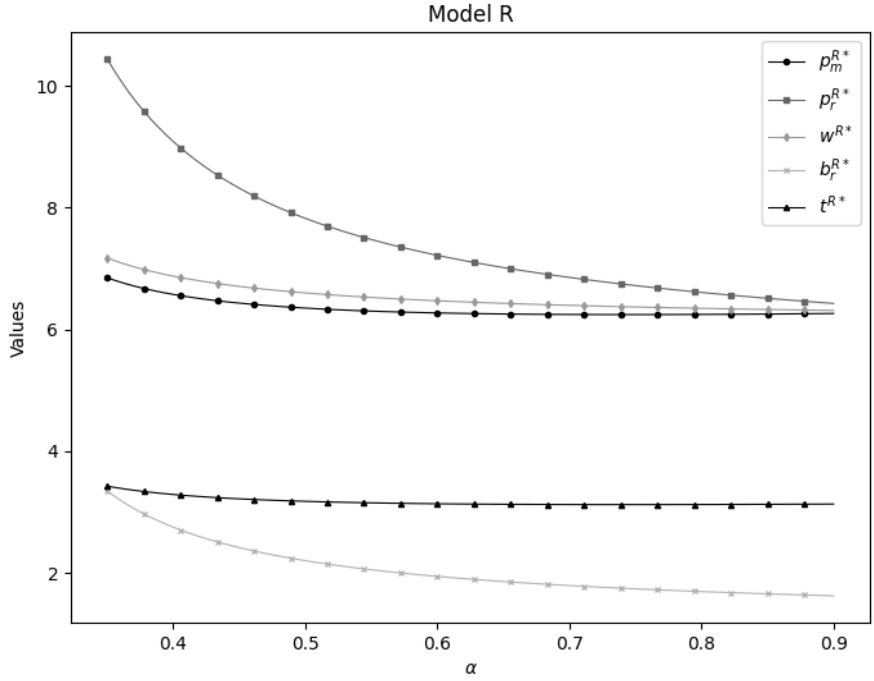

**FIGURE 4.** Decision Variables for Model R

Subsequently, **FIGURE 4** illustrates the decision variables for Model R, with the parameter $\alpha$ ranging from 0.35 to 0.9 to ensure stability in the calculations. The parameters employed are $s = 6$, $c_m = 10$, $c_r = 6$, and $\Delta = 4$. The decision variables depicted include $p_m^{R*}$, $p_r^{R*}$, $w^{R*}$, $b_r^{R*}$, and $t^{R*}$. The findings reveal a declining trend in most pricing variables as $\alpha$ increases, with a distinct pricing hierarchy of $p_r^{R*} > p_m^{R*}$. The subsidy $b_r^{R*}$ decreases as $\alpha$ rises, indicating a diminished reliance on retailer incentives in favor of direct sales. The transfer price $t_R^*$ also exhibits a downward trend, reflecting cost adjustments for used product transfers. Within a retailer-led recycling framework, these observations suggest that firms managing trade-in programs must refine their strategies to address an increase in direct-channel demand. Retailers might find it advantageous to enhance localized promotions or offer unique customer services rather than relying solely on substantial buyback subsidies. Such strategies could include value-added services like on-site assessment of used items or prioritized logistics for returned products, which can differentiate the retailer's offerings from the manufacturer's direct sales channel. Furthermore, renegotiating transfer prices with the manufacturer can help the retailer maintain adequate profit margins as incentives decrease. Maintaining a retail price above the manufacturer's direct-channel price is crucial to ensure the retailer's profitability, even when consumer preferences shift significantly toward the direct channel.

Finally, **FIGURE 5** delineates the decision variables for Model MR, with $\alpha$ ranging from 0.35 to 0.9 to ensure computational stability. The parameters are defined as $s = 6$, $c_m = 10$, $c_r = 6$, and $\Delta = 4$. The decision variables depicted include $p_m^{MR*}$, $w^{MR*}$, $b_m^{MR*}$, $b_r^{MR*}$, and $t^{MR*}$. The findings reveal a mixed trend, with certain variables increasing while others decrease as $\alpha$ varies. Notably, $p_m^{MR*}$ diverges, suggesting that remanufacturing-based hybrid models necessitate distinct pricing strategies. The subsidies $b_m^{MR*}$ and $b_r^{MR*}$ exhibit upward trends, underscoring the importance of financial incentives for the return of used products. The transfer price $t^{MR*}$ remains stable, indicating a consistent cost structure for product exchanges. Within this framework, it is crucial for both the manufacturer and the retailer to carefully coordinate subsidy adjustments to prevent margin erosion while aiming to increase the volume of trade-in items. For instance, if each party independently raises its trade-in subsidy, the total recycling rates may increase, but this could come at the expense of combined profitability. Instead, a collaboratively developed cost allocation strategy (e.g., assigning specific percentages to remanufacturing expenses or logistics fees) can facilitate the flow of more used items into remanufacturing process without causing excessive subsidy escalation. Consequently, within this mixed trend dynamic, it is essential to establish clear intervals for adjusting subsidies and to monitor cost thresholds, ensuring that cooperative relationships remain effective.

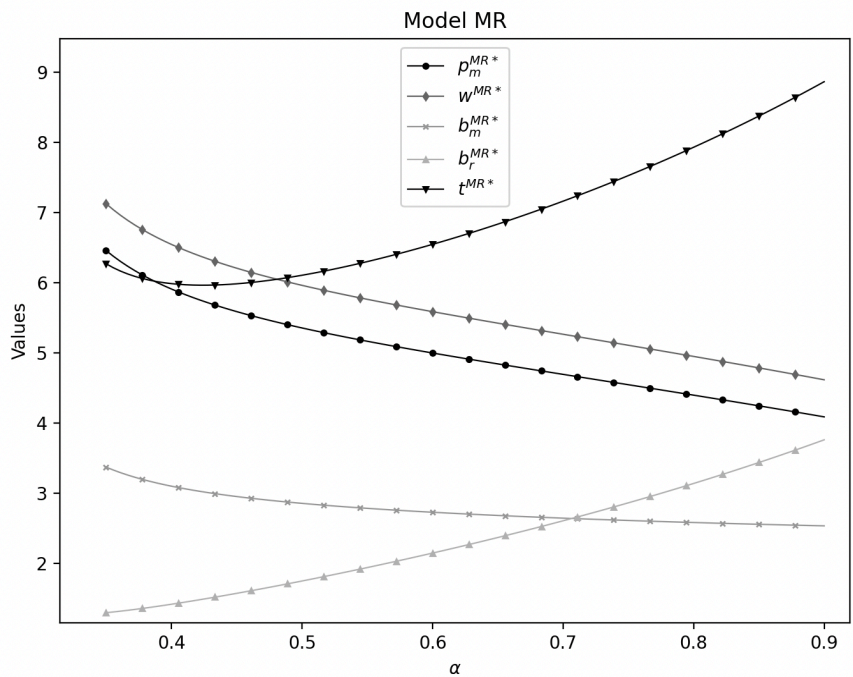

**FIGURE 5.** Decision Variables for Model MR

The managerial implications of these findings suggest that pricing and subsidy strategies should be dynamically adjusted in response to the growing customer preference for direct purchasing. As the inclination towards direct channels strengthens, manufacturers are advised to strategically reduce wholesale and transfer prices to maintain market share while optimizing profitability. Furthermore, in co-recycling framework, it is crucial for firms to meticulously balance pricing and subsidy frameworks across both direct and retail channels to maximize profitability and ensure cost-effective product return strategies. An adaptive pricing strategy that encompasses both direct and indirect channels will be vital for sustaining a competitive advantage in an evolving market environment.

# 5    Conclusions

This paper investigates optimal pricing strategies in DCCLSCs, explicitly addressing customer heterogeneity and comparing the implications of three distinct recycling frameworks: manufacturer-led, retailer-led, and collaborative recycling. By developing and applying a Manufacturer-Stackelberg game-theoretic model based on clearly defined customer demand functions for primary and replacement customers, we derive equilibrium prices and subsidies essential to effective supply chain management. Our model specifically distinguishes between primary customers, who purchase new products without engaging in returns, and replacement customers, who actively participate in trade-in programs by returning used products. By explicitly accounting for these differentiated customer segments, our model provides precise insights into pricing dynamics and recycling incentives within dual-channel structures.

Our analysis provides several novel findings that distinctly enhance existing research. Firstly, manufacturers consistently set lower prices in the direct sales channel compared to retailers, irrespective of the recycling structure. This strategic pricing decision is crucial, as it significantly influences consumer purchasing behavior, particularly among replacement customers who are highly sensitive to price differences and incentives associated with trade-in programs. By setting lower direct-channel prices, manufacturers effectively attract these replacement customers, thereby increasing product returns, which directly support both sustainability and profitability objectives. This finding contrasts with and significantly extends previous studies such as Chen and Wu (2024), who focused primarily on offline channel preferences and service level considerations without integrating the detailed dynamics between primary and replacement customer segments or the strategic implications of recycling incentives. By clearly highlighting these differences, our study fills a crucial research gap by providing a comprehensive, integrated analysis of pricing strategies influenced by customer heterogeneity.

Secondly, our analysis uncovers significant variations in the scale and necessity of trade-in subsidies across different recycling frameworks, an insight not thoroughly explored in prior literature. Manufacturer-led recycling emerges as a stable framework where moderate, targeted subsidies suffice to effectively stimulate product returns and maintain competitive pricing structures. In contrast, retailer-led recycling frameworks require considerably higher subsidies, reflecting the increased complexity and coordination challenges inherent in retailer-driven recycling initiatives. Such complexities are not adequately addressed in recent works by Yu et al. (2024, 2025), who primarily focused on centralized decision-making and platform-based sales modes without a nuanced examination of subsidy strategies across different recycling scenarios. Our research distinctly identifies collaborative recycling as superior, consistently yielding lower overall prices, significantly higher trade-in participation rates, and improved overall recycling efficiency. By explicitly comparing these frameworks, our research demonstrates clear managerial pathways for implementing effective recycling and pricing strategies tailored to specific operational contexts.

Thirdly, our research uniquely highlights the critical role of primary customers' preferences for direct-channel purchasing in shaping pricing and subsidy strategies. We find that stronger preferences for direct-channel purchases among primary customers correlate directly with reduced direct-channel prices and enhanced subsidies from manufacturers. This nuanced understanding significantly extends the insights provided by He et al. (2024) and Wei (2024), who studied competitive dynamics and subsidy impacts without explicitly integrating customer segmentation and channel preferences into their models. Our findings offer managers clear guidelines for optimizing pricing strategies to align with evolving consumer preferences, providing strategic leverage in highly competitive dual-channel environments.

From a managerial perspective, our results offer explicit, actionable strategies for contemporary businesses operating within dual-channel structures. Manufacturers should strategically manage direct-channel pricing and increase targeted subsidies in markets characterized by strong direct-channel preferences to effectively capture primary customers and maximize product returns. Retailers, particularly in retailer-led recycling scenarios, should concurrently adopt flexible pricing structures and carefully negotiate transfer pricing agreements with manufacturers to sustain profitability while encouraging robust trade-in participation. Our findings also provide clear and actionable policy recommendations. Policymakers should actively incentivize and facilitate collaborative recycling frameworks between manufacturers and retailers due to their demonstrated effectiveness in achieving balanced outcomes for profitability and sustainability. Real-world examples such as Apple's successful trade-in programs and Hewlett-Packard's comprehensive cartridge recycling initiatives illustrate the practical applicability of these insights, demonstrating how firms can operationalize theoretical recommendations to achieve significant environmental and economic benefits in real market conditions.

Despite its valuable contributions, this study acknowledges certain limitations, thereby opening avenues for future research. Our model currently assumes fully rational consumer behavior driven solely by economic incentives. Incorporating behavioral economic factors—such as consumer trust, perceived convenience, and emotional motivators—would enhance realism and applicability. Additionally, empirical validation through industry-specific case studies would greatly strengthen the practical relevance of our theoretical findings. Furthermore, examining multi-echelon closed-loop supply chains with additional stakeholders, including third-party recyclers and multiple distributors, could further enrich the strategic insights and operational implications derived from our model. In conclusion, by comprehensively integrating customer heterogeneity, direct-channel preferences, and multiple recycling frameworks into our pricing strategy analyses, this paper significantly advances the literature on DCCLSCs. Our findings provide practical managerial guidelines, explicit policy implications, and detailed strategic insights that empower businesses and policymakers to effectively navigate the complexities of dual-channel recycling, optimizing profitability and enhancing sustainability in contemporary commercial environments.